\documentclass[12pt]{article}
\usepackage{amsmath}    
\usepackage{amsfonts}
\usepackage{graphicx}   
\usepackage{verbatim}   
\usepackage{color}      
\usepackage{subfigure}  
\usepackage{hyperref}   

\begin{document}

\begin{center}
{\large \bf  Denis N. Sidorov and Nikolay A. Sidorov }\\
$ \,$\\
{\Large \bf Convex Majorants Method in the Theory of Nonlinear Volterra Equations~\footnote{This work partly funded by Russian Federal Framework Program ``Cadri''  P696, 30th of May 2010,
it is also supported by RFBR under the project \mbox{¹ 09--01--00377.} } }\\ %

\end{center}

\begin{abstract}
The main solutions  in sense of Kantorovich of nonlinear Volterra operator-integral equations are constructed. Convergence of the successive
approximations is established through studies of majorant integral and majorant algebraic equations. Estimates are derived for the
solutions and for the intervals on the right margin of which the solution has blow-up or solution start branching.
\end{abstract}



\section{Introduction}


Let us consider the following nonlinear continuos operator
$$\Phi (\omega_1, \dots, \omega_n,u,t): E_1\times \dots \times E_1 \times \mathbb{R}^1 \rightarrow E_2 $$
of $n+1$   variables  $\omega_1, \dots ,\omega_n, u, $  which are abstract continuos functions of real variable 
$t$ with values in $E_1.$ Here  $E_1, E_2$ are Banach spaces and $\Phi(0,\dots ,0, u_0, 0)=0,$
 $$K_i : \underbrace{\mathbb{R}^1 \times \dots  \mathbb{R}^1}_{i+1} \times \underbrace{E_1  \times \dots \times E_1}_{i} \rightarrow E_2$$ are nonlinear continuos operators
depending on the vector function $u(s) = (u(s_1),\dots , u(s_n))$  and  $t, s_1, \dots s_n$ are real variables.
Let $$\omega_i(t) = \int_0^t \cdots \int_0^t K_i(t,s_1,\dots , s_i, u(s_1), \dots , u(s_i)) ds_1 \dots ds_i,\, i=\overline{1,n}$$
   and address the following operator-integral equation for $t\in [0,T)$
   $$F(u,t) \equiv \Phi\biggl(  \int\limits_0^t K_1(t,s,u(s)) ds, \int\limits_0^t \int\limits_0^t K_2(t,s,s_1,s_2,u(s_1),u(s_2)) ds_1 ds_2, \dots $$ $$ \dots \int\limits_0^t \cdots \int\limits_0^t K_n(t,s_1, \dots s_n, u(s_1), \dots u(s_n)) ds_1 \dots ds_n,
   u(t), t  \biggr) = 0. \eqno{(1)}$$
Unknown abstract continuos function  $u(t)$ maps into  $E_1.$   Our objective is to find continuos solution $u(t) \rightarrow u_0$   from  $t \rightarrow 0.$
For $E_1 = E_2 = \mathbb{R}^1 $ the equation (1) has been studied by many authors in number of cases (see \cite{apartsyn2010, belbas}). 
However,  to the best of our knowledge, the equation (1) has not yet been studied  in general case of the Banach spaces $E_1, E_2$ . 

One of the common constructive methods in the theoretical and applied studies is method of majorants. 
Lev Kantorovich in his work \cite{kantor2} studied the functional equations in  $B_K$--spaces and converted the classical method of majorants into its 
abstract form, which makes it's methodology more clear and more unified. In his monograph (readers may refer to \cite{kantor}, p. 467) he specifically outlined the role
of the {\it main} solutions of nonlinear equations and role of the corresponding majorants. The main solutions are unique by definition and can be constructed 
using the successive approximations from the equivalent equation (4) starting from zero initial estimate. The main continuos solution $u^+(t)$ 
in the points of the interval $[0,T)$  satisfies  the estimate
 $|| u^+(t) ||_{E_1} \leq z(t)$ where $z(t)$ is continuos positive solution of the following majorant Volterra integral equation
 $$z(t)=f\biggl( \int_0^t \gamma(z(s)) ds\biggr), \,\, z(t) \in C_{[0,T)}^+. \eqno{(2)}$$
 Here and below $f, \gamma$ are monotone increasing continuos functions.
 If the trivial solution $u=0$ satisfies the equation (1) then this solution is the main solution.
 This trivial case is excluded from the consideration below. 
 If one continue the main nontrivial solution $u^+(t)$ outside of the the interval $[0,T)$ (where we see the convergence of the successive approximations) in the right hand side direction from the margin 
 point $T$  then the solution $u^+(t)$ can go to $\infty$ or start {\it branching} [5].
 Obviously, there is a case when operator $F$ satisfies the Lipschitz condition for $\forall u $ and the main solution is continuable on the whole interval  $[0,\infty).$
 If the Lipschitz condition is not fulfilled, then in addition  to the main solution the equation (1) can have arbitrarily  many other continuous solutions which
 cross the main solution.\\
 
  \noindent {\it Example 1}.\\
    $$ u(t) = p\int_0^t u^{\frac{p}{p-1}}(s), \, 1<p<\infty$$
    Here $u_1(t)=0$  is the main solution. Other continuous solutions:    $u_2(t)=t^p,$
    $$ u_c(t) = \left\{ \begin{array}{ll}
         0, & \mbox{ $-\infty \leq t \leq c$}\\
        (t-c)^p, & \mbox{ $c\leq t <\infty$}\end{array} \right. $$
  The main solution $u^+(t)=0$ of this example is singular solution for the correspondent Cauchy problem $\dot{u}=pu^{\frac{p}{p-1}}, \, u(0)=0.$
         
The objective of this paper is to construct main solutions for the equation (1) on the maximal interval $[0,T).$ The paper consists of two parts,  illustrative examples and conclusion. 
In the second part for the equation (1) existence theorem is derived for the main solution $u(t) \rightarrow u_0$ when $t \rightarrow 0$ with  estimate  $|| u(t) ||_{E_1}$  when $t\in [0,T).$

 We propose the approach for construction of the approximations  $u_n(t)$   and the interval  $[0,T)$ on which they converge point-wise for  $\forall u_0(t), $
 if $||u_0(t)||_E \leq z^+(t),$  where $z^+(t)$ is the main nonnegative solution of the corresponding majorant integral equation (2).
 The sufficient conditions are derived if $\lim\limits_{t\rightarrow T} z^+(t) = \infty$
 (or $\lim\limits_{t\rightarrow T} \frac{dz^+(t)}{dt} = +\infty), $   i.e. the main solution of the majorant equation (or its derivative) has
 the  blow-up limit (goes to $\infty$ for finite time T).
 Under such a conditions the unknown solution $u(t)$ of the equation (1) can also strives to infinity during the finite time $T^{\prime} \geq T$
 or appear to be branching.
  
  In the third part of this paper we demonstrate how to construct and employ the following
  majorant algebraic systems
      $$ \left\{ \begin{array}{ll}
         \mbox{$r = R(r,t)$}, \\
         \mbox{$1 = R_r^{\prime} (r,t)$}\end{array} \right. \eqno{(3)} $$
         for construction of the main solution of the equation (1).
         In the algebraic system (3) $R(0,0)=0, \, R_r^{\prime}(0,0)=0, \, R(r,t)$ is the convex function wrt $r.$
         The algebraic majorant systems (3)  were also called as the Lyapunov majorants~\cite{greben}.
         Such majorants as well as more general algebraic majorants were used in mechanics (see [6], p.198--216)  and for the construction
         of implicit functions in spaces $B_K$ .
         It is to be noted that algebraic majorant systems has the unique positive solution  $r^*, T^*.$
         
          Using this approach one can define the {\it guaranteed interval} $[0,T^*],$  on which the equation (1) has
         the main solution $u(t) \rightarrow 0$ for $t\rightarrow 0$ and radius  of the sphere  $S(0,r^*)$   in the space  $C_{[0,T^*]}^{E_1},$    
         in which the main solution can be constructed using the successive approximations which converge uniformly.
          
          \section{Integral majorants in construction of\\ the main solution}
          
         Let's consider the following equation 
          $$u=\mathcal{L}(u), \eqno{(4)}$$
        ãäå $ \mathcal{L}(u) = A^{-1} (Au - F(u,t)), $
        which is equivalent to the equation (1). Here 
  $A$ is continuously invertible  operator from   $E_1$ into $E_2.$     If the operator $F$ has the Frechet derivative  $F_u(0,0)$
   and its invertible then we can assume   $A=F_u(0,0).$
         
 \noindent {\bf Definition.}\\  If the approximations $u_n(t)=\mathcal{L}(u_{n-1}), \, u_0=0$
for $t\in [0,T^*)$  strives to the solution $u^+(t)$ of the equation (4), then function  $u^+(t)$  we call Kantorovich main solution of the equation (1).

It is to noted that here we follow the monograph of Lev V. Kantorovich (refer to \cite{kantor}, p.467), where the term   {\it ``the main solution of the functional equation''}
has been formulated~\cite{kantor, kantor2}. Under the solution we will be assuming the main solution below in this paper.

Let us study the operator  $F(u,t)-Au:$   $C_{[0,T]}^{E_1} \rightarrow C_{[0,T]}^{E_2}.$
Here $C_{[0,T]}^{E_1}$ and $C_{[0,T]}^{E_2}$ are complete spaces.

We will get  the estimate in norms of the spaces  $E_1, E_2$:\\
{\bf A)}  $||F(u,t) - Au||_{E_2} \leq f(\int_0^t \gamma(||u(s)||_{E_1}))ds), \, t \in [0,T).$\\
Let in the inequality  {\bf A)} and below the following assumption be hold:\\
{\bf B)} $\gamma, \, f$  are continuously -- monotone  increasing functions on the the segments  $[0,z^{\prime}]$  and  $[\gamma(0), \gamma(z^{\prime})],\, z\leq \infty$
correspondingly;\\
{\bf C)} for $t \in [0,T)$ exist the function $z^{\prime}(t)$ in the cone  $C_{[0,T]}^+$  such as
 $$ z^{\prime}(t) \geq f\biggl( \int_0^t \gamma(z^{\prime}(s))ds \biggr). \eqno{(5)}$$
     
 \noindent {\bf Remark 1.}\\  In the Lemmas 2 and 3 we propose the method to define the margin  $T$ such as for  $ t \in [0,T)$
 the condition {\bf C)}  will be fulfilled.
 Because of the condition {\bf A)} $f(\gamma(0)t) \geq 0.$
 Zero is lower solution of the majorant integral equation (2), and $z^{\prime}(t)$ is upper solution in the cone $C_{[0,T)}^+.$

Under the conditions {\bf C)} and {\bf B)} we introduce the sequence  $$ z_n(t) = f\biggl( \int_0^t \gamma(z_{n-1}(s))ds \biggr), \, z_0=0. $$
 Then due to the Theorem  2.11 (see \cite{kantor}, p.464) for   $\forall n, \, t\in [0,T)$ the inequalities 
  $$0=z_0(t)\leq z_1(t) \leq \cdots \leq z_n(t) \leq z^{\prime}(t)$$ are fulfilled.
  Hence the limit   $\lim\limits_{n\rightarrow \infty} z_n(t) = z^+(t)$ exist.
  Since  $\gamma, \, f$ are continuos functions and due to the Lebeg theorem  (see, e.g. , \cite{shilov})  the limit exist
  $$ \lim_{n\rightarrow \infty} f \biggl( \int_0^t \gamma(z_n(s))ds \biggr) = f\biggl( \int_0^t \gamma(z^{+}(s))ds \biggr). $$  
  Thus function$z^+(t)$ appears to be continuos on  $[0,T)$  and to be the main solution of the majorant equation (2).
 Approximation $z_n(t)$ in the points of the interval $[0,T)$ converge to $z^+(t), \, z^+(t) \in C_{[0,T]}^+.$
 
 Let us now proceed to the construction of the solution $u^+(t)$ of the equation (1) using the successive approximations.
 In addition to {\bf A)}, {\bf B)} è {\bf C)} let the following inequality be fulfilled\\
 {\bf D)} $||F(u+\Delta u,t) - F(u,t) - Au ||_{E_2} \leq f \biggl( \int_0^t \gamma(||u(s)||_{E_1} + ||\Delta u(s)||_{E_1})ds \biggr) - $\\ $- f(\int_0^t \gamma (||u(s)||_{E_1})ds).$\\
 Under the condition of the Frechet differentiability of the operators  $$F(u,t), \, f\biggr(\int\limits_0^t \gamma(z(s)) ds \biggl)$$ the verification of the 
 inequality  {\bf D)}  can be replaced with  verification of the condition {\bf E)} (see below).
 Indeed, let mentioned Frechet derivatives exist and continuos for  $t\in [0,T)$ following the norms of linear bounded operators in the spaces
  $\mathcal{L}(E_1 \rightarrow E_2)$  and  $\mathcal{L}(C_ {[0,T)}^+ \rightarrow C_{[0,T)}^+)$ correspondingly.
   Under such assumption we assume functions $f, \, \gamma$  has monotone increasing  and continuos derivatives and the Frechet  differential   $f$ is defined by the formula
     $$f_z^{\prime} \biggl( \int\limits_0^t \gamma (z(s) ds  \biggr) h \equiv f_{\gamma}^{\prime} \biggl( \int\limits_0^t \gamma (z(s) ds  \biggr)  \int\limits_0^t \gamma_z^{\prime} (z(s)) h(s) ds  $$
     for   $\forall h(s) \in C_{[0,T]}^+.$
     
     Let in addition to the conditions {\bf A)} and {\bf B)}  $\forall V(t) \in C_{[o,T]}^{E_1}$ the following inequality be fulfilled\\
     {\bf E)} $||(F_u(u,t) - A)V ||_{E_2} \leq f_z^{\prime} \biggl( \int\limits_0^t \gamma (||u(s)||_{E_1} ds  \biggr) ||V||_{E_1},\,$ $||V||_{E_1} \in C_{[0,T)}^+.$
     Then we have the following lemma  
    
    \noindent {\bf Lemma 1.}\\
Let the inequality {\bf E)} be fulfilled and the derivatives  $f_{\gamma}^{\prime}, \, \gamma_z^{\prime}$ are  monotone increasing.  Then inequality  {\bf D)} be fulfilled.

{\it Proof.}  Let us employ the Lagrange finite-increments formula (\cite{vat}, p.367) and conditions of the Lemma 1. Then we get the inequality 
$$||F(u+\Delta u,t) - F(u,t) - A \Delta u ||_{E_2} =$$ $$= ||\int_0^1 (F_u (u+\Theta \Delta u, t) - A) d \Theta \Delta u||_{E_2} \leq \int_0^1 f_{\gamma}^{\prime} \biggl ( \int_0^t \gamma ( ||u(s)||_{E_1}+$$
$$+ \Theta ||\Delta u (s)||_{E_1} )ds \biggr) \int_0^t \gamma^{\prime} \biggl (  ||u(s)||_{E_1} +\Theta ||\Delta u(s) ||_{E_1} \biggr) ||\Delta u(s)||_{E_1} ds d\Theta = $$
$$= f(\int_0^t \gamma (||u(s)||_{E_1} + ||\Delta u (s)||_{E_1}) ds ) - f (\int_0^t \gamma (||u(s)||_{E_1})ds).$$

Let us now construct the approximations $u_n(t) = \mathcal{L}(u_{n-1}), \, u_0 = 0$
to solution  $u^+(t).$   We follow the proof  of the Theorem 2.22 (\cite{kantor}, p.466) 
  and state the estimates
  $||u_{n+p}(t) - u_n(t)||_{E_1} \leq z_{n+p}(t) - z_n(t)$ 
   for $t\in [0,T),$ where $z_n(t) = f (\int_0^t \gamma (z_{n-1}(s)) ds), \, z_n(t) \in C_{[0,T)}^+, \, u_n(t)\in C_{[0,T)}^{E_1}, \, u_0=0, z_0=0.$
   
   It is to be noted here that similar estimates in different problem has been also used in our paper  \cite{sid_kluwer} for studies of explicit
   mappings based on convex majorants method.
   
   Due to the conditions  {\bf A)}, {\bf B)} and {\bf C)} and  following the above mentioned approach the limit $\lim\limits_{n\rightarrow \infty} z_n(t) = z^+(t)$
   exists for 
   $\forall t \in [0,T),$ i. å.  $z_n(t)$ is fundamental sequence in the each point   $t\in [0,T).$
    Hence the sequence of abstract functions  $u_n(t)$  with values in the Banach space  $E_1$ for each 
    $t \in [0,T)$ converges in norms of the space   $E_1$  to function $u^+(t).$  Because of  the operator  $\mathcal{L}(u)$
    is continuos then the equality $ u^+(t) = \mathcal{L}(u^+)$ is fulfilled,
    i.e.  $u^+(t)$ satisfies the equation (1) and belongs to the space $C_{[0,T)}^{E_1}.$
     
   Hence the following theorem be fulfilled    

\noindent {\bf Theorem 1.}\\  Let for   $t\in [0,T)$ the conditions {\bf A)}, {\bf B)}, {\bf C)} and {\bf D)} are fulfilled.
Then equation  (1) in the space $C_{[0,T)}^{E_1}$   has main solution $u^+(t).$
Moreover,  $||u^+(t)||_{E_1} \leq z^+(t),$ were $z^+(t)$ is main solution of majorant equation (2), approximations $u_n(t)=\mathcal{L}(u_{n-1}), \, u_0=0$
converge to  $u^+(t)$ in norm of the space  $E_1$  for  $\forall t \in [0,T^+),$  approximations  $z_n(t) = f(\int_0^t \gamma (z_{n-1}(s))ds ), z_0=0$
converge to  $z ^+(t).$

  In Theorem 1 $T^+$  remains not defined. For the definition of  $T^+$ we reduce the majorant integral equation (2) to Cauchy theorem for 
separable differential equation. For this objective we introduce the differentiable function 
  $\omega(t) = \int\limits_0^t \gamma (z(s)) ds.$   Òîãäà $\frac{d \omega(t)}{d t} = \gamma (z(t)), \, \omega(0)=0,$
  where  $z(t) = f(\omega(t)).$
  That is why the Cauchy problem which is equivalent to the equation (2) is following 
        $$ \left\{ \begin{array}{ll}
         \mbox{$\frac{d \omega}{ dt } = \gamma (f(\omega (t)))$}, \\
         \mbox{$\omega(0) = 0.$}\end{array} \right. \eqno{(2^{\prime})} $$
  Lemma 2  and Lemma 3 defines the estimate of the interval  $[0,T^+),$
on which the Cauchy problem  ($2^{\prime}$) in space $C_{[0,T^{+})}^+$
has the unique solution  $\omega^+(t)$ and approximations   $\omega_n(t) = \int\limits_0^t \gamma (f(\omega_{n-1}(s)))ds, \, \omega_0=0$
converge to this unique solution.

\noindent {\bf Lemma 2.}\\
Let $\gamma (f(\omega))$ be continuos, strictly positive and monotone increasing function. Let exists
$\lim\limits_{\omega \rightarrow \infty} \int\limits_0^{\omega} \frac{d \omega}{ \gamma (f(\omega))} = T^+.$
Then ($2^{\prime}$) in cone $C_{[0,T^+]}^+$
has monotone increasing solution $\omega^+(t).$ The approximations $\omega_n(t) = \int\limits_0^t \omega (f(\omega_{n-1}(s))) ds, \, \omega_0=0$
converge to  $\omega^+(t),$  $\lim\limits_{t \rightarrow T^+} \omega^+ (t) = \infty.$

\noindent  {\it Proof}.  Let us separate the variables in  ($2^{\prime}$) and reduce the Cauchy problem to search for the positive monotone 
increasing branch of implicit function 
$\omega = \omega(t), \, \omega(0)=0$ from equation  $\Phi(\omega)=t,$  where $\Phi(\omega) = \int\limits_0^{\omega} \frac{d \omega}{\gamma(f(\omega))}.$
If $\gamma (f(\omega))$ is rational fraction then  antiderivative $\Phi(\omega)$ can be explicitly constructed in terms of logarithms, arctangences and  rational functions.
 It is to be noted that under conditions of the lemma 1, function $\Phi(\omega)$ is continuos and monotone increasing on semi-axis $[0,\infty),$
  i.e. $\Phi^{\prime} = \frac{1}{\gamma(f(\omega))} > 0,$
   $\lim\limits_{\omega \rightarrow 0} \Phi(\omega) = 0, \, \lim\limits_{\omega \rightarrow \infty} \Phi(\omega) = T^{+}.$
   Hence the  mapping  $\Phi:  [0,\infty) \rightarrow [0,T^+)$
   is bijective, equation $\Phi (\omega) = t$ for $0\leq t < T^+$
   uniquely defines function $\omega^+(t),$ which obviously  satisfies the  integral equation 
   $$\omega(t) = \int\limits_0^t \gamma (f(\omega(s)))ds. $$
   Because of the monotone increasing of the functions  $f$ and $\gamma,$ the approximations $\omega_n(t) = \int_0^t \gamma (f(\omega_{n-1}(s)))ds, \, \omega_0 = 0$
   for $t \in [0,T^+)$ converge to $\omega^+(t).$

If  $\gamma(f(\omega))$ is rational fraction, then in number of cases the solution $\omega^+(t)$ can be explicitly constructed in complicated cases using the computer 
algebra systems \cite{apartsyn2010}.

\noindent   {\bf Remark 2.}  \\
  For known $\omega^+(t)$ using the formula  $z^+(t)=f(\omega^+(t))$ we find the solution of majorant integral equation (2).
  It is to be note that under conditions of the Lemma 2, approximations  $z_n= f(\int\limits_0^t \gamma(z_{n-1}(s)ds)), \, z_0=0$
converge for  $t\in [0,T^+)$ to the solution $z^+(t)$.

\noindent {\bf Remark 3.} \\
 If under the conditions of Lemma 2 then $\lim\limits_{\omega \rightarrow \infty} \int\limits_0^{\omega} \frac{d\omega}{\gamma(f(\omega))} = \infty$
then the solution $z^+(t)$ is continuable on $[0, \infty).$
This result follows from the Theorem 2.7 (\cite{barb}, p.~148).

For example, let inequality 
$$||F(u,t)-Au||_{E_2} \leq a \int\limits_0^t ||u(s)||ds + b, \, a>0, \, b>0,$$
be fulfilled for  $\forall u, \, 0\leq t < \infty.$
Then majorant integral equation  (2) will be linear as follows  
$z(t)=a \int\limits_0^t z(s) ds +b$
and has the unique solution $z(t)=b e^{at}, \, 0\leq t < \infty.$
Çàìåòèì, ÷òî â ýòîì ñëó÷àå $\gamma (f(\omega)) = a\omega + b$
$\lim\limits_{\omega \rightarrow \infty} \int\limits_0^{\omega} \frac{d\omega}{a\omega + b} = \infty.$
If in this case
 $$||F(u+\Delta u, t) - F(u,t) - A\Delta u ||_{E_2} \leq a \int\limits_a^t ||\Delta u(s)||_{E_1} ds, $$
  then conditions of the theorem 1 are fulfilled on semi-axis $0\leq t < \infty$
  and equation (1) will have the solution $u^+(t)$ in the space $C_{[0,\infty)}^{E_1},$
  $|| u^+(t)||_{E_1} \leq b e^{at}.$  Obviously, from this result not follows the fact that in area $||u(t)||_{E_1} \geq b e^{at}$ 
  the equation (1) does not have another solutions.
  
 \noindent {\bf Lemma 3.}\\
 Let superposition $\gamma (f(\omega))$ be continuos and strictly positive for  $0\leq \omega \leq \omega^*.$
 Let limits $\lim\limits_{\omega \rightarrow \omega^*} \gamma (f(\omega)) = \infty $ are exists.
 $\lim\limits_{\omega \rightarrow \omega^*} \int\limits_0^{\omega} \frac{d \omega}{\gamma (f(\omega))} = T^+.$
 Then Cauchy problem  ($2^{\prime}$) for $t\in [0,T^+]$
 in the cone $C_{[0,T^+]}^+$ has continuos monotone increasing solution $\omega^+(t),$
 and $\lim\limits_{t\rightarrow T^+} \frac{d\omega^+}{dt} = 0,$
 approximations $\omega_n(t) = \int_0^t \gamma (f \omega_{n-1} (s)) ds, \, \omega_0=0$
 converge for $0\leq t \leq T^+$
 to the solution $\omega^+(t).$

   Proof of the lemma 3 follows from the bijectivity of mapping $\Phi:\, [0,\omega^*] \rightarrow [0,\Phi(\omega^*)]$
for  $\Phi(\omega) = \int_0^{\omega} \frac{d\omega}{\gamma(f(\omega))},\,\, \Phi(\omega^*) = T^+. $

\noindent  {\bf Remark 4.}\\ Under the conditions of Lemma 3 the point $T^+$ is  blow-up limit of the derivative of solution $z^+(t)$
  of majorant equation (2).

  \section{Algebraic majorants in construction of\\ the main solution}

 Let in the equation (1) $u_0=0,$ I.e. $\Phi(0,\dots,0)=0.$  
Our objective is to construct continuous solution  $u^+(t)$
with successive approximations $u_n(t) = \mathcal{L}(u_{n-1})$ 
in close interval $[0,T^+].$ In the space $C_{[0,T^+]}^{E_1}$
we introduce the norm $||u||=\max\limits_{0\leq t\leq T^+} ||u(t)||_{E_1}.$
We suppose, that operator $F$ is Frechet differentiable w.r.t. $u$
Let for $0\leq t \leq T,$ ãäå $T\geq T^+$ è $u\in S(0,r) \subset E_1, $ inequalities be fulfilled :\\
$\mathbf{A^{\prime})}$  $||F(u,t)-Au||_{E_2} \leq f(r,t);$\\
$\mathbf{E^{\prime})}$ $||F^{\prime}_u(u,t) - A||_{E_2} \leq f_r^{\prime}(r,t);$\\
$\mathbf{G)}$ 
Let functions $f(r,t), \, f_r^{\prime}(r,t)$ are positive $r>0, \, t>0$
and monotone increase, $f(0,0)=0,\, f_r^{\prime}(0,0) \in [0,1),$
function $f(r,t)$ convex w.r.t. $r $.
Then algebraic equation $r=||A^{-1}|| f(r,t)$ according to definition 5.1 from monograph [3], p.~205.
will be the Lyapunov majorant for operator  $\mathcal{L}(u).$ Because of monotone increasing function 
$f(r,t), \, f_r(r,t)$ and convexity of the function  $f(r,t)$ system 
  $$ \left\{ \begin{array}{ll}
       \mbox{$r = ||A^{-1}|| f(r,t)$}, \\
       \mbox{$1 = ||A^{-1}|| f_r^{\prime}(r,t)$}\end{array} \right. $$
        has unique positive solution $r^+, T^+.$
        Moreover, equation $r=||A^{-1}|| f(r,t)$ where $0\leq t \leq T^+$ ([3], p. 218) uniquely defined monotone increasing solution $r^*=r(t)$
        Approximations $r_n(t) = ||A^{-1}|| f(r_{n-1}(t),t), \, r_0=0,$ when $0\leq t \leq T^+$ converge to the function $r(t).$  Corresponding approximations $r_n = ||A^{-1}|| f(r_{n-1},T^+), \, r_0=0,$
converge to $r^+.$  Function  $r(t)$  is main solution of the  Lyapunov  majorante equation.
        On the base of the Lemma 5.1 ([3],  page 206)
        if  $||u_i(t)||_{E_1} \leq r_i, \, i=1,2,$ $||u_2(t) - u_1(t)|| \leq r_2 - r_1,$ then for $0\leq t \leq T^+$ 
        $$ || \mathcal{L}(u_2) - \mathcal{L}(u_1)||_{E_1} \leq ||A^{-1}|| (f(r_2,t)-f(r_1,t)).$$
         Apart from approximations $r_n(t)$ for  solutions for Lyapunov majorante, we introduce approximation $u_n(t)=\mathcal{L}(u_{n-1}), \, u_0=0$
         of the main solution of the equation (1).
          For arbitrary $k$ and $l\geq k$ because of the conditions $\mathbf{A^{\prime})}, \, \mathbf{E^{\prime})}$  and above mentioned inequality, we come to estimate $||u_l(t) - u_k(t)||_{E_1} \leq r_l(t) - r_k(t) \leq r_l(T^+) - r_k(T^+).$
          Such that  $r_l(T^+)$ monotone increasing sequence and $\lim\limits_{l\rightarrow \infty} r_l(T^+) = r^+,$
          then $||u_l(t)-u_k(t)||_{E_1} \leq \varepsilon$ for $l,k \geq N(\varepsilon)$ if  $t \in [0,T^+].$
          Hence  $||u_l(t)-u_k(t)||_{E_1} \leq \varepsilon$
          for   $l,k \geq N(\varepsilon).$  Because of complete space $C_{[0,T^+]}^{E_1}$ exist limit 
          $\lim\limits_{l\rightarrow \infty} u_l(t) = u^+(t).$
          Moreover, $u^+(t)$ continuos w.r.t. $t,$ and  approximation $u_n(t) = \mathcal{L}(u_{n-1}), \, u_0=0$
     converge on segment $[0,T^+]$  uniformly w.r.t. $t.$
          
          Then follows \\
\noindent {\bf Theorem 2.} \\
Let  $\Phi(0,\dots, 0)=0,$  inequalities $\mathbf{A^{\prime})}$, $\mathbf{E^{\prime})}$  are satisfied  when $t \in [0,T^+], $ pair  $(r^+, T^+), \, r^+>0,\, T^+>0$
satisfies algebraic system 
  $$ \left\{ \begin{array}{ll}
       \mbox{$r = ||A^{-1}|| f(r,t)$}, \\
       \mbox{$1 = ||A^{-1}|| f_r^{\prime}(r,t)$,} \end{array} \right.$$
       where function $f(r,t)$ satisfies the condition $\mathbf{G)}$  
Then on $[0,T^+]$ equation (1) has continuous solution $u^+(t)$ in space $C_{[0,T^+]}^{E_1}.$
More over, approximations $u_n(t) = \mathcal{L}(u_{n-1})$
converge uniformly w.r.t. $t,$ $\max\limits_{0\leq t \leq T^+} || u^+(t)|| \leq r^+.$

\noindent {\it Example 2.}\\ 
Let us consider the following problem
    $$ \left\{ \begin{array}{ll}
       \mbox{$\frac{\partial^2 u(x,t)}{\partial x^2} + \int_0^t \sin (t-\tau + x) u^2(x,\tau) d\tau = t$}, \\
       \mbox{$u\bigl|_{x=0} = u\bigl|_{x=1} =0, \,\,\, 0\leq x \leq 1, \, t\geq 0.$} \end{array} \right.$$
       We search for classical solution $u \rightarrow 0$ for $t\rightarrow 0.$
        Here $E_1 = 
       \hspace{-.3cm}\stackrel{\hspace*{.4cm}\circ\, (2)}{C}_{\hspace*{-.3cm}[0,1]}$
       --  space of twice differentiable w.r.t.   $x$ functions are zero on the margins  $[0,1],\, E_2=C_{[0,1]}.$ $Au = \frac{\partial^2 u}{\partial x^2,}$ operator $A \in \mathcal{L}(E_1\rightarrow E_2)$
        has limited reverse $A^{-1} = \int\limits_0^1 G(x,s) [\cdot] ds,$
        where  $$G(x,s) = \left\{ \begin{array}{ll}
       x(s-1), & \mbox{$0\leq x\leq s \leq 1$}, \\
       s(x-1), & \mbox{$s\leq x \leq 1$,} \end{array} \right.$$
$|| A^{-1}||_{\mathcal{L}(E_1\rightarrow E_2)} \leq 1.$

 Following the Theorem 1, the corresponding majorant integral equation (2) is following  $z(t) = \int_0^t z^2(s) ds + t.$
Then function $z^+(t) = \tan t\,$ for $ 0\leq t < \frac{\pi}{2}$ is the main 
solution of majorant integral equation. Therefore $\frac{\pi}{2}$
 is point , in which has blow-up solution limit $z^+(t).$
 Boundary problem on the base of Theorem 1 has in space 
$C_{[0,\frac{\pi}{2})}^{E_1}$ the solution $u^+(x,t),$  beside  $0\leq t < \frac{\pi}{2}$
   $$ \max\limits_{0\leq x \leq 1}  \biggl ( \biggl| \frac{\partial^i u^+(x,t) }{\partial x^i}\biggr |, i=0,1,2 \biggr ) \leq {\tan} \, t.$$
From other hand, if we follow the Theorem 1, we construct majorant algebraic equation $r = tr^2 + t.$
Function $r^+(t) = \frac{1-\sqrt{1-4t^2}}{2t}$ for $t \in[0,0.5]$ is the main solution of majorant algebraic equation.

    According to Theorem 2 we construct the following system 
         $$ \left\{ \begin{array}{ll}
       \mbox{$r = t r^2 + t$}\\
       \mbox{$1=2 t r $} \end{array} \right.$$
  which has one positive solution $T^+=0.5, \, r^+=1.$
    Therefore according to the theorem 2 we get guaranteed interval  w.r.t.  $t$ of existence of the solution $u^+(x,t)$
    of the boundary problem $[0,1/2]$ with estimate of the norm of the solution $u^+$ such as 
   $$\max\limits_{x\in[0,1],\,\, t\in[0,0.5]}  \biggl\{ \biggl | \frac{\partial^i u^+(x,t)}{\partial x^i} \biggr|, i=0,1,2  \biggr\} \leq 1,\, ||u^+||_{E_1} \leq r^+(t), \, 0 \leq t \leq 0.5.$$
    Since $0.5 <\frac{pi}{2}$ then in this example integral majorant provide more precise estimate $u^+,$ comparing to the algebraic  one.

    As the footnote let us notice that with studies of the equation (1) in $B_K$  spaces and with introduction of the abstract norms in the Kantorovich sense 
    it is possible to get more fine systems of majorant integral and algebraic equations. Such a majorants will  characterize the solution of equation (1)
    more deeply.
    Majorant algebraic equations possible to construct and it is possible to study the solutions of $n$-diminutional Volterra equations (1), namely for $t\in \mathbb{R}^n, \, n\geq 2.$  
    As a matter of fact, the algebraic majorants provides more  rough estimates comparing to the integral majorants, but its more easy to constract and to employ the algebraic majorants.
    Since the solution of the majorant integral equation has blow-up limit, for numerical solution in the neighborhood of such points it make sense to employ
    the adaptive meshes.

\vfill

\hfill 

\end{document}